\documentclass{article}

\usepackage{amsmath, amssymb, amsthm, amsfonts, bm}
\usepackage{enumerate}

\numberwithin{equation}{section}
\newtheorem{theorem}{Theorem}[section]
\newtheorem{cthm}{Theorem}
\newenvironment{ctheorem}[1]
  {\renewcommand\thecthm{#1}\cthm}
  {\endcthm}

\newtheorem{ccor}{Corollary}
\newenvironment{ccorollary}[1]
  {\renewcommand\theccor{#1}\ccor}
  {\endccor}
\newtheorem{lemma}{Lemma}[section]
\theoremstyle{definition}
\newtheorem{definition}{Definition}[section]
\theoremstyle{remark}

\newcommand{\Rnum}[1]{\uppercase\expandafter{\romannumeral #1\relax}}
\newcommand{\rnum}[1]{\romannumeral #1\relax}
\newcommand{\mr}[1]{\mathrm{#1}}
\newcommand{\mb}[1]{\mathbb{#1}}
\newcommand{\mc}[1]{\mathcal{#1}}

\DeclareMathOperator{\av}{Av}
\DeclareMathOperator{\mt}{\mathcal{M}_\mathcal{T}}
\DeclareMathOperator{\md}{\mathcal{M}_d}

\title{A sharp integral inequality for the dyadic maximal operator and a related stability result}
\author{Eleftherios N. Nikolidakis}
\date{}

\begin{document}
\maketitle
\footnotetext{ {\em E-mail address}: enikolid@uoi.gr}

\begin{abstract}
We prove a sharp integral inequality for the dyadic maximal operator due to which the evaluation of the Bellman function of this operator with respect to two variables is possible, as can be seen in \cite{3}. Our inequality of interest is proved in this article by a simpler and more immediate way. We also study a stability result in connection with this inequality, that is we provide a necessary and sufficient condition, for a sequence of functions, under which we obtain equality in the limit. The proof of this result is based on the proof of the related inequality which we present in this article.
\end{abstract}

\section{Introduction} \label{sec:1}
The dyadic maximal operator on $\mb R^n$ is a useful tool in analysis and is defined by
\begin{equation} \label{eq:1p1}
\md \phi(x) = \sup \left\{ \frac{1}{|Q|}\int_Q|\phi(y)|\,\mr dy: x\in Q,\ Q\subseteq \mb R^n\ \text{is a dyadic cube}\right\},
\end{equation}
for every $\phi\in L_\text{loc}^1(\mb R^n)$, where the dyadic cubes are those formed by the grids $2^{-N}\mb Z^n$, for $N=0, 1, 2,\ldots$.
As is well known it satisfies the following weak type (1,1) inequality
\begin{equation} \label{eq:1p2}
\left|\left\{x\in\mb R^n : \md \phi(x) > \lambda\right\}\right| \leq
\frac{1}{\lambda} \int_{\{\md\phi>\lambda\}}|\phi(y)|\,\mr dy,
\end{equation}
for every $\phi\in L^1(\mb R^n)$ and every $\lambda>0$, from which it is easy to get the following $L^p$-inequality
\begin{equation} \label{eq:1p3}
\|\md\phi\|_p \leq \frac{p}{p-1}\|\phi\|_p,
\end{equation}
for every $p>1$ and $\phi\in L^p(\mb R^n)$.

It is easy to see that the weak type inequality \eqref{eq:1p2} is best possible. It has also been proved that \eqref{eq:1p3} is best possible (see \cite{1}, \cite{2} for general martingales and \cite{18} for dyadic ones).

For the study of the dyadic maximal operator it is desirable for one to find refinements of the above mentioned inequalities. Concerning \eqref{eq:1p2}, improvements have been given in \cite{10} and \cite{11}. If we consider \eqref{eq:1p3}, there is a refinement of it if one fixes the $L^1$-norm of $\phi$. That is we wish to find explicitly the following function (named as Bellman) of two variables $f$ and $F$.
\begin{equation} \label{eq:1p4}
B_Q^{(p)}(f,F) = \sup\left\{ \frac{1}{|Q|} \int_Q (\md \phi)^p : \phi \geq 0,\ \frac{1}{|Q|}\int_Q\phi=f,\ \frac{1}{|Q|}\int_Q\phi^p=F\right\},
\end{equation}
where $Q$ is a fixed dyadic cube and $f, F$ are such that $0< f^p\leq F$.

This function was first evaluated in \cite{5}. In fact it has been explicitly computed in a much more general setting of a non-atomic probability space $(X,\mu)$ equipped with a tree $\mc T$, with structure  similar to the one that the dyadic subcubes of $[0,1]^n$ have (see the definition in Section \ref{sec:2}). Then we define the associated maximal operator by
\begin{equation} \label{eq:1p5}
\mt\phi(x) = \sup\left\{ \frac{1}{\mu(I)}\int_I |\phi|\,\mr d\mu: x\in I\in \mc T\right\},
\end{equation}
for every $\phi\in L^1(X,\mu)$.

Moreover, \eqref{eq:1p2} and \eqref{eq:1p3} still hold in this setting and remain sharp. Now if we wish to refine \eqref{eq:1p3} we should introduce the so-called Bellman function of the dyadic maximal operator of two variables given by
\begin{equation} \label{eq:1p6}
B_{\mc T}^{(p)}(f,F) = \sup\left\{ \int_X(\mt\phi)^p\,\mr d\mu: \phi\geq 0,\ \int_X\phi\,\mr d\mu = f,\ \int_X\phi^p\,\mr d\mu = F\right\},
\end{equation}
where $0< f^p\leq F$.
This function of course generalizes \eqref{eq:1p4}. In \cite{5} it is proved that
\[
B_{\mc T}^{(p)}(f,F) = F\,\omega_p\!\left(\frac{f^p}{F}\right)^p,
\]
where $\omega_p : [0,1] \to \bigl[1,\frac{p}{p-1}\bigr]$ is defined by $\omega_p(z) = H_p^{-1}(z)$, and $H_p(z)$ is given by $H_p(z) = -(p-1)z^p + pz^{p-1}$. As a consequence $B_{\mc T}^{(p)}(f,F)$ does not depend on the structure of the tree $\mc T$.
The technique for the evaluation of \eqref{eq:1p6}, which is used in \cite{5}, is based on an effective linearization of the dyadic maximal operator that holds on an adequate class of functions called $\mc T$-good (see the definition in Section 2), which is enough to describe the problem that is settled on \eqref{eq:1p6}.
In \cite{8} now a different approach has been given for the evaluation of \eqref{eq:1p6}. This was actually done for the Bellman function of three variables in a different way avoiding the calculus arguments that are given in \cite{4}. More precisely the following is a consequence of the results in \cite{8}.

\begin{ctheorem}{A} \label{thm:a}
Let $\phi\in L^p(X,\mu)$ be non-negative with $\int_X\phi\,\mr d\mu=f$. Then the following inequality is true and sharp
\begin{equation} \label{eq:1p7}
\int_X(\mt\phi)^p\,\mr d\mu \leq -\frac{1}{p-1}f^p + \frac{p}{p-1}\int_X \phi\,(\mt\phi)^{p-1}\,\mr d\mu.
\end{equation}
\end{ctheorem}
This inequality as one can see in \cite{8}, enables us to find a direct proof for the exact evaluation of \eqref{eq:1p6}. For this evaluation we also need a symmetrization principle that can be found in \cite{8} (presented as Theorem \ref{thm:2p1} below) and which is also used in this article for the sharpness of our results.  In this paper we will prove the following generalization of Theorem \ref{thm:a}.
By using the linearization technique that appears in \cite{5} in a more complicated form, we present in Section 3 a proof of the theorem that appears just below (mentioned as Theorem 1), which generalizes Theorem A and which is the following.
\begin{ctheorem}{1} \label{thm:1}
Let $\phi$ be as in the hypothesis of Theorem A and suppose that $q\in [1,p]$. Then the following inequality is true for any $\beta>0$
\begin{multline} \label{eq:1p8}
\begin{aligned}&\int_X(\mt\phi)^p\,\mr d\mu \leq -\frac{q(\beta+1)}{(p-1)q\beta+(p-q)}f^p + \end{aligned} \\
+\frac{p(\beta+1)^q}{(p-1)q\beta+(p-q)}\int_X\phi^q(\mt\phi)^{p-q}\,\mr d\mu.
\end{multline}
Additionally \eqref{eq:1p8} is best possible for any given $q\in[1,p]$, $f>0$ and $\beta$ such that $0<\beta \leq \frac{1}{p-1}$. By this we mean that if one fixes the second constant appearing on the right hand side of inequality \eqref{eq:1p8} then we cannot increase the absolute value of the first constant appearing in front of $f^p$ in a way such that \eqref{eq:1p8} still holds.
\end{ctheorem}
The following is also true and is an easy consequence of Theorem 1.

\begin{ccorollary}{1} \label{cor:1}
Let $\phi: (X,\mu)\to \mb R^+$ be such that $\int_X\phi\,\mr d\mu=f$. Then for every $q\in [1,p]$ the following inequality holds
\begin{equation} \label{eq:1p9}
\int_X(\mt\phi)^p\,\mr d\mu \leq -\frac{q}{p-1}f^p + \left(\frac{p}{p-1}\right)^q\int_X\phi^q(\mt\phi)^{p-q}\,\mr d\mu.
\end{equation}
Additionally \eqref{eq:1p9} is best possible for any given $q\in[1,p]$ and $f>0$.
\end{ccorollary}

Moreover, by using the symmetrization principle that is mentioned below (Theorem 2.1) and Theorem 1 we easily derive inequalities of Hardy type as described by the following
\begin{ccorollary}{2} \label{cor:2}
For any $g: (0,1]\to\mb R^+$ non-increasing such that $\int_0^1 g(u)\,\mr du=f$, the following inequality is true for any $\beta>0$ and sharp for any
$\beta$ such that $0<\beta \leq\frac{1}{p-1}$.
\begin{multline} \label{eq:1p10}
\begin{aligned} &\int_0^1\left(\frac{1}{t}\int_0^tg(u)\,\mr du\right)^p\mr dt \leq
-\frac{q(\beta+1)}{(p-1)q\beta+(p-q)}f^p +  \end{aligned}\\
+\frac{p(\beta+1)^q}{(p-1)q\beta+(p-q)}\int_0^1\left(\frac{1}{t}\int_0^t g(u)\,\mr du\right)^{p-q}\!\!g^q(t)\,\mr dt.
\end{multline}
\end{ccorollary}
For the case $q=1$ and the value $\beta=\frac{1}{p-1}$ inequality \eqref{eq:1p10} is well known and is in fact equality, as can be seen by applying a simple integration by parts argument. We also note that inequality \eqref{eq:1p8} is also a consequence of the results in \cite{3}, where it is proved a more general inequality which involves also the parameter $A=\int_X\phi^q\,\mr d\mu$. In this paper we ignore this parameter and give a more direct proof of \eqref{eq:1p8}.

Moreover the proof that we give for \eqref{eq:1p8} enables us to provide a stability result for this inequality. That is we characterize when we do have equality in the limit in \eqref{eq:1p8} for a sequence of functions $(\phi_n)_n$. More precisely we prove the following
\begin{ctheorem}{2} \label{thm:2}
Let $(\phi_n)_n$ be a sequence of nonnegative, $\mc T$-good functions (the exact definition will be given in Section 2) satisfying $\int_X\phi_n\,\mr d\mu=f$ and $\int_X\phi_n^p\,\mr d\mu=F$, for every $n\in N$ and $q$ be such that $q\in(1,p)$. Let also $\beta$ which satisfies
\begin{equation} \label{eq:1p11}
\beta+1= \omega_p\!\left(\frac{f^p}{F}\right)
\end{equation}
Then $(\phi_n)_n$ satisfies equality in the limit in \eqref{eq:1p8}, if and only if the following is true
\begin{equation} \label{eq:1p12}
\underset{n}{\lim}\int_X|\mt\phi_n-(\beta+1)\phi_n|^p\,\mr d\mu=0,
\end{equation}
\end{ctheorem}

By using now the results of \cite{9} we conclude that if we fix the $L^1$ and $L^p$ norms of $\phi_n$, $n\in N$, the sequence $(\phi_n)_n$ gives equality in the limit in \eqref{eq:1p8} for any $q\in(1,)]$, if and only if it behaves as an extremal sequence for the respective Bellman function \eqref{eq:1p6}.
At last we mention that the evaluation of \eqref{eq:1p6} has been given by an alternative method in \cite{13} while certain Bellman functions
corresponding to several problems in harmonic analysis, have been studied in \cite{6}, \cite{7}, \cite{14}, \cite{15}, \cite{16} and \cite{17}.

\section{Preliminaries} \label{sec:2}
Let $(X,\mu)$ be a non-atomic probability space. We give the following from \cite{5} or \cite{8}.

\begin{definition} \label{def:2p1}
{\em A set $\mc T$ of measurable subsets of $X$ will be called a tree if the following are satisfied}
\begin{enumerate}[i)]{\em
\item $X\in\mc T$ and for every $I\in\mc T$, $\mu(I) > 0$.
\item For every $I\in\mc T$ there corresponds a finite or countable subset $C(I)$ of $\mc T$ containing at least two elements such that
\vspace{-5pt}}
\begin{enumerate}[a)]{\em
\item the elements of $C(I)$ are pairwise disjoint subsets of $I$
\item $I = \bigcup\, C(I)$.}
\end{enumerate}{\em
\item $\mc T = \bigcup_{m\geq 0} \mc T_{(m)}$, where $\mc T_{(0)} = \left\{ X \right\}$ and
\[
\mc T_{(m+1)} = \bigcup_{I\in \mc T_{(m)}} C(I).
\]
\item The following holds
\[
\lim_{m\to\infty} \sup_{I\in \mc T_{(m)}} \mu(I) = 0
\]
\item The tree $\mc T$ differentiates $L^1(X,\mu)$.
}
\end{enumerate}
\end{definition}
The last property stated in the definition of the tree $\mc T$ means that for every $\phi\in L^1(X,\mu)$,
$\lim_{x\in I\in \mc T, \mu(I)\to0}\frac{1}{\mu(I)}\int_I\phi\,\mr d\mu=0$ for $\mu$-almost all $x$ on $X$.

For the proof of Theorem \ref{thm:1} we will use an effective linearization for the operator $\mt$ that was introduced in \cite{5}. We describe it as appears there and use it in the sequel.

For every $\phi\in L^1(X,\mu)$  non-negative and $I\in\mc T$ we define $\av_I(\phi) = \frac{1}{\mu(I)}\int_I\phi\,\mr d\mu$. In the proofs below we denote
$\av_I(\phi)$ by $y_I$.
We will say that $\phi$ is $\mc T$-good if the set
\[
\mc A_\phi = \left\{x\in X: \mt \phi(x) > \av_I(\phi)\ \text{for all}\ I\in\mc T\ \text{such that}\ x\in I\right\}
\]
has $\mu$-measure zero.

Let now $\phi$ be $\mc T$-good and $x\in X\!\setminus\!\mc A_\phi$. We define $I_\phi(x)$ to be the largest in the nonempty set
\[
\left\{ I\in\mc T: x\in I\ \text{and}\ \mt\phi(x) = \av_I(\phi)\right\}.
\]
Now given $I\in\mc T$ let
\begin{align*}
A(\phi,I) &= \left\{x\in X\!\setminus\!\mc A_\phi: I_{\phi}(x) = I\right\}\subseteq I\ \ \text{and} \\
S_\phi &= \left\{I\in\mc T: \mu(A(\phi,I))>0\right\} \cup \left\{X\right\}.
\end{align*}
Obviously then
\[
\mt\phi = \sum_{I\in S_\phi} \av_I(\phi) \chi_{A(\phi,I)},\ \mu\text{-a.e.},
\]
where $\chi_E$ is the characteristic function of $E$.
We also define the following correspondence $I\to I^\star$ by: $I^\star$ is the smallest element of $\left\{J\in S_\phi: I\subsetneq J\right\}$. It is defined for every $I\in S_\phi$ except $X$. Also it is obvious that the $A(\phi,I)$'s are pairwise disjoint and that
\[
\mu\left(\bigcup_{I\notin S_\phi} A(\phi,I)\right) = 0,
\]
so that
\[
\bigcup_{I\in S_\phi} A(\phi,I) \approx X,
\]
where by $A\approx B$ we mean that
\[
\mu(A\!\setminus\!B) = \mu(B\!\setminus\!A) = 0.
\]
Now the following is true (see \cite{5}).
\begin{lemma} \label{lem:2p1}
Let $\phi$ be $\mc T$-good
\begin{enumerate}[i)]
\item If $I, J\in S_\phi$ then either $A(\phi,J)\cap I = \emptyset$ or $J\subseteq I$.
\item If $I\in S_\phi$ then there exists $J\in C(I)$ such that $J\notin S_\phi$.
\item For every $I\in S_\phi$ we have that $I\approx \underset{\substack{J\in S_\phi\\J\subseteq I}}{\bigcup} A(\phi,J)$.
\item For every $I\in S_\phi$ we have that
\[
A(\phi,I) = I \setminus \underset{\substack{J\in S_\phi\\ J^\star=I}}{\bigcup} J,
\]
so that
\[
\mu(A(\phi,I)) = \mu(I) - \sum_{\substack{J\in S_\phi\\ J^\star=I}} \mu(J).
\]
\end{enumerate}
\end{lemma}
\noindent From the above we see that
\[
\av_I(\phi) = \frac{1}{\mu(I)} \sum_{\substack{J\in S_\phi\\ J\subseteq I}} \int_{A(\phi,J)}\phi\,\mr d\mu.
\]
In the sequel we will also need the notion of the decreasing rearrangement of a $\mu$-measurable function defined on $X$. This is given by the following equation
\[
\phi^\star(t) = \sup_{\substack{e\subseteq X\\ \mu(e) \geq t}} \Bigl[ \inf_{x\in e}|\phi(x)| \Bigr],\ \ t\in (0,1].
\]
This is the unique non-increasing, left continuous function defined on $(0,1]$, equimeasurable to $|\phi|$ (that is $\mu\!\left(\left\{|\phi|>\lambda\right\}\right) = \left|\left\{\phi^\star>\lambda\right\}\right|$, for any $\lambda>0$). A more intuitive definition of $\phi^\star$ is that it describes a rearrangement of the values of $|\phi|$ in decreasing order.
We are now ready to state the following theorem, which appears in \cite{8}, and can be viewed as a symmetrization principle for the dyadic maximal operator.
\begin{theorem} \label{thm:2p1}
The following equality is true
\begin{multline} \label{eq:2p1}
\begin{aligned} &\sup\left\{ \int_K G_1(\mt\phi)\, G_2(\phi)\,\mr d\mu: \phi^\star = g,\ \phi\geq 0,\right. \\
&\hspace{30pt} \left. \vphantom{\int_K} K\, \text{measurable subset of}\ X\ \text{with}\ \mu(K)=k\right\} = \end{aligned} \\
=\int_0^k G_1\!\!\left(\frac{1}{t}\int_0^t g\right) G_2(g(t))\,\mr dt,
\end{multline}
where $G_i: [0,+\infty) \to [0,+\infty)$ are increasing functions for $i=1,2$, while $g: (0,1]\to \mb R^+$ is non-increasing. Additionally the supremum in \eqref{eq:2p1} is attained by some $(\phi_n)$ such that $\phi_n^\star = g$, for every $n\in N$. This sequence of functions is independent of the pair of functions $(G_1, G_2)$.
\end{theorem}

\section{Proof of the inequality \eqref{eq:1p8}} \label{sec:3}

We now proceed to the
\begin{proof}[Proof of Theorem \ref{thm:1}] ~ \\
Let $\phi:(X,\mu)\to\mb R^+$ be $\mc T$-good such that $\int_X\phi\,\mr d\mu=f$ and let $q\in (1,p]$. (The case $q=1$ can be handled easily if we consider a sequence $(q_n)_n$ of elements of $(1,p]$, tending to $q=1$ and applying the result for every $q_n$). We consider the quantity
\[
k_q = \int_X\phi^q (\mt\phi)^{p-q}\,\mr d\mu.
\]
By the definition of the linearization of the dyadic maximal operator we have that
\begin{equation} \label{eq:3p1}
k_q = \sum_{I\in S_\phi} \int_{A(\phi,I)} \phi^q\,\mr d\mu \cdot y_I^{p-q}.
\end{equation}
By H\"{o}lder's inequality now, since $q>1$, we have that
\begin{equation} \label{eq:3p2}
\int_{A(\phi,I)}\phi^q\,\mr d\mu \geq \frac{1}{\alpha_I^{q-1}}\biggl(\int_{A(\phi,I)}\phi\,\mr d\mu\bigg)^q,
\end{equation}
where $A(\phi,I) = I\!\setminus\!\bigcup_{J\in S_\phi,J^\star=I}J$, in view of Lemma \ref{lem:2p1} \rnum 4), and so $\alpha_I = \mu(A(\phi,I)) = \mu(I)-\sum_{J\in S_\phi, J^\star=I}\mu(I)$. Thus \eqref{eq:3p1} in view of \eqref{eq:3p2} gives
\begin{align} \label{eq:3p3}
k_q &\geq \sum_{I\in S_\phi} y_I^{p-q} \frac{\left(\int_I\phi\,\mr d\mu - \sum_{J\in S_\phi, J^\star=I}\int_J\phi\,\mr d\mu\right)^q}{\left(\mu(I) - \sum_{J\in S_\phi, J^\star=I}\mu(J)\right)^{q-1}} = \notag \\
&= \sum_{I\in S_\phi}y_I^{p-q} \frac{\left(\mu(I)y_I - \sum_{J\in S_\phi, J^\star=I} \mu(J)y_J\right)^q}{\left(\mu(I) - \sum_{J\in S_\phi, J^\star=I}\mu(J)\right)^{q-1}}.
\end{align}
We use now H\"{o}lder's inequality in the following form
\begin{equation} \label{eq:3p4}
\frac{(\lambda_1 + \lambda_2 + \ldots + \lambda_m)^q}{(\sigma_1 + \sigma_2 + \ldots + \sigma_m)^{q-1}} \leq \frac{\lambda_1^q}{\sigma_1^{q-1}} + \frac{\lambda_2^q}{\sigma_2^{q-1}} + \ldots + \frac{\lambda_m^q}{\sigma_m^{q-1}},
\end{equation}
which holds for every $\lambda_i\geq 0,\ \sigma_i>0$ since $q>1$. \\
We consider now an arbitrary $\beta$ such that $0\leq\beta\leq\frac{1}{p-1}$ . We set for any $I\in S_\phi$
\[
\tau_I = (\beta+1)-\beta\rho_I,\ \ \text{where}\ \ \rho_I = \frac{\mu(A(\phi,I))}{\mu(I)} = \frac{\alpha_I}{\mu(I)},
\]
thus concluding that $\tau_I>0$. For this choice of $\tau_I$ we have that
\begin{equation} \label{eq:3p5}
\tau_I \mu(I) \ -\  (\beta+1)\sum_{J\in S_\phi, J^\star=I}\mu(J) = \mu(I) \ - \sum_{J\in S_\phi, J^\star=I}\mu(J).
\end{equation}
Thus using \eqref{eq:3p4} and \eqref{eq:3p5} we have from \eqref{eq:3p3} that
\begin{align} \label{eq:3p6}
k_q &\geq \sum_{I\in S_\phi} y_I^{p-q} \Biggl\{ \frac{(\mu(I) y_I)^q}{(\mu(I)\tau_I)^{q-1}} - \sum_{\substack{J\in S_\phi\\ J^\star=I}} \frac{(\mu(J)y_J)^q}{((\beta+1)\mu(J))^{q-1}}\Biggr\} = \notag \\
 &= \sum_{I\in S_\phi}\mu(I)\frac{y_I^p}{\tau_I^{q-1}} - \sum_{I\in S_\phi}y_I^{p-q}\sum_{\substack{J\in S_\phi\\ J^\star=I}} \frac{y_J^q}{(\beta+1)^{q-1}}\mu(J).
\end{align}
By the definitions now of $S_\phi$ and the correspondence $I\to I^\star$ for $I\neq X$, we conclude from \eqref{eq:3p6} that
\begin{align} \label{eq:3p7}
k_q &\geq \sum_{I\in S_\phi}\mu(I)\frac{y_I^p}{\tau_I^{q-1}} - \sum_{\substack{I\in S_\phi\\ I\neq X}} \frac{1}{(\beta+1)^{q-1}}y_I^q(y_{I^\star})^{p-q}\mu(I) = \notag \\
 &= \sum_{I\in S_\phi}\frac{1}{\rho_I}\alpha_I\frac{y_I^p}{((\beta+1)-\beta \rho_I)^{q-1}} - \frac{1}{p}\sum_{\substack{I\in S_\phi\\ I\neq X}}\frac{py_I^q(y_{I^\star})^{p-q}}{(\beta+1)^{q-1}}\mu(I).
\end{align}
We now use the following elementary inequality
\[
p x^q\!\cdot\! y^{p-q} \leq q x^p + (p-q) y^p,
\]
which holds since $1< q\leq p$ for any $x, y>0$. By \eqref{eq:3p7} we thus have
\begin{align} \label{eq:3p8}
k_q &\geq \sum_{I\in S_\phi} \frac{\alpha_I}{\rho_I} \frac{y_I^p}{((\beta+1)-\beta \rho_I)^{q-1}} - \frac{1}{p}\sum_{\substack{I\in S_\phi\\ I\neq X}}\frac{\left[qy_I^p + (p-q)(y_{I^\star})^p\right]}{(\beta+1)^{q-1}} \mu(I) = \notag \\
&= \sum_{I\in S_\phi}\frac{\alpha_I}{\rho_I} \frac{y_I^p}{((\beta+1)-\beta \rho_I)^{q-1}} - \frac{p-q}{p}\frac{1}{(\beta+1)^{q-1}} \sum_{\substack{I\in S_\phi\\ I\neq X}} (y_{I^\star})^p\mu(I) - \notag \\
&\hspace{100pt} -\frac{q}{p}\frac{1}{(\beta+1)^{q-1}}\sum_{I\in S_\phi}y_I^p\mu(I) + \frac{q}{p}\frac{1}{(\beta+1)^{q-1}}y_X^p.
\end{align}
By using now Lemma 2.1 iv), and the definition of the correspodence $I\to I^\star$, we have that
\[
\sum_{\substack{I\in S_\phi\\ I\neq X}}(y_{I^\star})^p\mu(I) = \sum_{I\in S_\phi}y_I^p(\mu(I)-\alpha_I),
\]
thus \eqref{eq:3p8} gives
\begin{multline*}
k_q \geq \sum_{I\in S_\phi}\frac{\alpha_I}{\rho_I}\frac{1}{((\beta+1)-\beta \rho_I)^{q-1}}y_I^p - \frac{p-q}{p}\frac{1}{(\beta+1)^{q-1}}\sum_{I\in S_\phi}(\mu(I)-\alpha_I)y_I^p - \\
- \frac{q}{p}\frac{1}{(\beta+1)^{q-1}}\sum_{I\in S_\phi}\mu(I)y_I^p + \frac{q}{p}\frac{1}{(\beta+1)^{q-1}}y_X^p.
\end{multline*}
After some simple cancellations we conclude that
\begin{multline} \label{eq:3p9}
k_q \geq \sum_{I\in S_\phi}\frac{\alpha_I}{\rho_I} \left(\frac{1}{((\beta+1)-\beta \rho_I)^{q-1}} - \frac{1}{(\beta+1)^{q-1}}\right)y_I^p + \\
+ \frac{p-q}{p}\frac{1}{(\beta+1)^{q-1}}\sum_{I\in S_\phi}\alpha_Iy_I^p + \frac{q}{p}\frac{1}{(\beta+1)^{q-1}}y_X^p.
\end{multline}
Now note that
\[
\frac{1}{((\beta+1)-\beta x)^{q-1}} - \frac{1}{(\beta+1)^{q-1}} \geq \frac{(q-1)\beta x}{(\beta+1)^q},
\]
by the mean value theorem on derivatives for all $x\in [0,1]$, so by \eqref{eq:3p9} we have as a consequence that
\begin{align} \label{eq:3p10}
k_q \geq& \sum_{I\in S_\phi} \left[\frac{\alpha_I}{\rho_I}\frac{(q-1)\beta \rho_I}{(\beta+1)^q}\right]y_I^p + \frac{p-q}{p}\frac{1}{(\beta+1)^{q-1}}\sum_{I\in S_\phi}\alpha_Iy_I^p + \frac{q}{p}\frac{1}{(\beta+1)^{q-1}}y_X^p \notag \\
 =& \sum_{I\in S_\phi}\left[\frac{(q-1)\beta}{(\beta+1)^q} + \frac{p-q}{p}\frac{1}{(\beta+1)^{q-1}}\right]\alpha_Iy_I^p + \frac{q}{p}\frac{1}{(\beta+1)^{q-1}}f^p,
\end{align}
and we have derived inequality \eqref{eq:1p8} for $\mc T$-good functions.

For the general $\phi:(X,\mu)\to\mb R^+$ which belongs to $L^p(X,\mu)$ we argue as follows. Consider the sequence $(\phi_m)_m$ defined by $\phi_m= \sum_{\substack{I\in \mc T_{(m)}}} \av_I(\phi)\chi_I$, and for any $m\in N$ set
$$\Phi_m= \sum_{\substack{I\in \mc T_{(m)}}} max\{\av_J(\phi):I\subseteq J\in \mc T\}\chi_I=\mt\phi_m.$$
The last equality holds due to the fact that $\av_J(\phi_m)=\av_I(\phi_m)=\av_I(\phi)$ whenever $J\subseteq I\in \mc T_{(m)}$. It is easy to see that $\int_X\phi_m\,\mr d\mu=\int_X\phi\,\mr d\mu=f$ while
$\int_X\phi^p_m\,\mr d\mu\leq\int_X\phi^p\,\mr d\mu$ for all $m$ and that $\Phi_m$ increases to $\mt\phi$ on $X$. Now $\phi_m$
satisfies \eqref{eq:1p8} since as can be easily seen is $\mc T$-good, and since $\mc T$ differentiates $L^1(X,\mu)$ we get that $\phi_m$ tends almost everywhere to $\phi$. Thus by taking limits and using the dominated convergence theorem we obtain \eqref{eq:1p8} for $\phi$.

At this point we give the following.
\begin{proof}[Proof of Corollary \ref{cor:2}]~\\
Let $g: (0,1]\to \mb R^+$ be non-increasing such that $\int_0^1g(u)\,\mr du=f$. Fix a non-atomic probability space $(X,\mu)$ equipped with a tree structure $\mc T$ that differentiates $L^1(X,\mu)$. Applying Theorem \ref{thm:2p1} for the pair of functions
\[
\left(G_1(t) = t^p, G_2(t)=1\right)\quad \text{and}\quad \left(G_3(t) = t^{p-q}, G_4(t)=t^q\right)
\]
we conclude that there exists, for every $n\in N$, $\phi_n: (X,\mu)\to\mb R^+$ such that $\phi_n^\star=g$, such that
\begin{equation} \label{eq:3p11}
\lim_n\int_X(\mt\phi_n)^p\,\mr d\mu = \int_0^1\left(\frac{1}{t}\int_0^tg\right)^p\mr dt,
\end{equation}
and
\begin{equation} \label{eq:3p12}
\lim_n\int_X\phi_n^q(\mt\phi_n)^{p-q}\,\mr d\mu = \int_0^1\left(\frac{1}{t}\int_0^tg\right)^{p-q}g^q(t)\,\mr dt.
\end{equation}
Applying \eqref{eq:1p8} for every $(\phi_n)$ and taking the limits as $n\to\infty$, we conclude by \eqref{eq:3p11} and \eqref{eq:3p12} the validity of inequality \eqref{eq:1p10}.

We now prove that \eqref{eq:1p10} is best possible. We proceed to this as follows: We first treat the case where $\beta=\frac{1}{p-1}$. We consider the following continuous and decreasing function $g_\alpha(t) = c\, t^{-\alpha}$, defined in $(0,1]$, where $c = f(1-\alpha)$, and
$\alpha\in \bigl(0,\frac{1}{p}\bigr)$.
Then it is easy to show that $\int_0^1 g_\alpha(u)\,\mr du=f$ while $g_\alpha\in L^p((0,1])$.

Note that for any $t\in (0,1]$ the equality $\frac{1}{t}\int_0^tg(u)\,\mr du = (\frac{1}{1-\alpha})g(t)$ is true. So considering the difference
\[
J = \int_0^1\left(\frac{1}{t}\int_0^tg_\alpha\right)^p\mr dt - \left(\frac{p}{p-1}\right)^q\int_0^1g_\alpha^q(t)\left(\frac{1}{t}\int_0^tg_\alpha\right)^{p-q}\mr dt
\]
we see that it is equal to
\[
J = \left(\frac{1}{1-\alpha}\right)^p\int_0^1g_\alpha^p(t)\,\mr dt - \left(\frac{p}{p-1}\right)^q\left(\frac{1}{1-\alpha}\right)^{p-q}\int_0^1g_\alpha^p(t)\,\mr dt.
\]
Since $\int_0^1g_\alpha^p(t)\,\mr dt = f^p(1-\alpha)^p\frac{1}{1-\alpha p}$ we have that
\begin{align*}
J &= \frac{f^p}{1-\alpha p} - \left(\frac{p}{p-1}\right)^q(1-\alpha)^q\frac{f^p}{1-\alpha p} = \\
&=-\frac{f^p}{1-\alpha p}\left[\left(\frac{p}{p-1}\right)^q(1-\alpha)^q-1\right] = -f^p \,G(\alpha),
\end{align*}
where $G(\alpha)$ is defined for any $\alpha\in \bigl(0,\frac{1}{p}\bigr)$ by $G(\alpha) = \frac{\left(\frac{p}{p-1}\right)^q(1-\alpha)^q-1}{1-\alpha p}$.
But as it is easily seen by using de L' Hospital's rule,
\[
\lim_{\alpha\to 1/p^-}G(\alpha) = -q\left(1-\frac{1}{p}\right)^{q-1}\left(\frac{p}{p-1}\right)^q\left(-\frac{1}{p}\right) = \frac{q}{p-1}.
\]
We now prove the sharpness of \eqref{eq:1p10} for any $\beta$ such that $0<\beta<\frac{1}{p-1}$.

We fix such a $\beta$ and we consider the following continuous and decreasing function $g_\beta(t) = c\, t^{-\alpha}$, defined on $(0,1]$, where $c = f(1-\alpha)$, and $\alpha=\frac{\beta}{\beta+1}$. Then
$\alpha\in \bigl(0,\frac{1}{p}\bigr)$,
and it is easy to see that $\int_0^1 g_\beta(u)\,\mr du=f$ while for any $\beta$ as above we have that $g_\beta\in L^p((0,1])$.

Moreover, $\int_0^1g^p_\beta(u)\,\mr du=\frac{f^p}{(\beta+1)^p} \frac{\beta+1}{1-\beta(p-1)}$. Note that for any $t\in (0,1]$ the following equality holds $\frac{1}{t}\int_0^tg_\beta(u)\,\mr du = (\beta+1)g_\beta(t)$. We then consider the difference
\begin{multline}
\begin{aligned}&J =\int_0^1g_\beta^q(t)\left(\frac{1}{t}\int_0^tg_\beta\right)^{p-q}\mr dt- \end{aligned} \\
-\left[\frac{(q-1)\beta}{(\beta+1)^q}+\frac{p-q}{p}\left(\frac{1}{\beta+1}\right)^{q-1}\right]\int_0^1\left(\frac{1}{t}\int_0^tg_\beta\right)^p\mr dt
\end{multline}
Then we easily see after some simple calculations that
$J=\frac{q}{p}\frac{1}{(\beta+1)^{q-1}}f^p$. The proof of Corollary 2 is now complete.
\end{proof}
Now for the proof of Theorem \ref{thm:2} we need to prove the sharpness of \eqref{eq:1p8}. This is easy now to show since by Theorem \ref{thm:2p1} for any $g: (0,1]\to \mb R^+$ non increasing, there exists a sequence $\phi_n: (X,\mu)\to\mb R^+$ of rearrangements of $g$ such that
\begin{equation} \label{eq:3p13}
\lim_n\int_X(\mt\phi_n)^p\,\mr d\mu = \int_0^1\left(\frac{1}{t}\int_0^tg\right)^p\mr dt
\end{equation}
and
\begin{equation} \label{eq:3p14}
\lim_n\int_X\phi_n^q (\mt\phi_n)^{p-q}\,\mr d\mu = \int_0^1g^q(t)\left(\frac{1}{t}\int_0^tg\right)^{p-q}\mr dt.
\end{equation}
We discuss now the case where $0<\beta<\frac{1}{p-1}$ and we consider the function $g_\beta$ (denoted now as $g$), constructed in the proof of Corollary \ref{cor:2}. For every $n\in N$, we choose a rearrangement $\phi_{n}$ of $g$ such that
\[
\left|\int_0^1\left(\frac{1}{t}\int_0^tg\right)^p\mr dt - \int_X\left(\mt\phi_{n}\right)^p\mr d\mu\right| \leq \frac{1}{n}
\]
and
\[
\left|\int_0^1g^q(t)\left(\frac{1}{t}\int_0^tg\right)^{p-q}\mr dt - \int_X\phi_{n}^q\left(\mt\phi_{n}\right)^{p-q}\mr d\mu\right| \leq \frac{1}{n}
\]
Then by the choice of $g$, we conclude that \eqref{eq:1p8} is best possible.
The case $\beta=\frac{1}{p-1}$ is entirely similar so we omit it.
The proof of Theorem \ref{thm:1} is now complete.
\end{proof}

\section{Proof of Theorem \ref{thm:2}} \label{sec:4}
\begin{proof}
We begin by stating inequality \eqref{eq:1p8} in the following equivalent form
\begin{equation} \label{eq:4p1}
\int_X(\mt\phi)^{p-q}\phi^q \,\mr d\mu\geq A_0(\beta)\int_X(\mt\phi)^p\,\mr d\mu+\frac{q}{p}\frac{1}{(\beta+1)^{q-1}}f^p,
\end{equation}
where $A_0(\beta)=\frac{(q-1)\beta}{(\beta+1)^q}+\frac{p-q}{p}\frac{1}{(\beta+1)^{q-1}}$.
For the one direction of the proof we suppose that $(\phi_n)_n$ satisfies $\int_X\phi_n\,\mr d\mu=f$ and $\int_X\phi_n^p\,\mr d\mu=F$, for every $n\in N$. Since \eqref{eq:1p11} is equivalent to
\begin{equation} \label{eq:4p2}
 F(\beta+1)^{p-q}-A_0(\beta)F(\beta+1)^{p}=\frac{q}{p}\frac{1}{(\beta+1)^{q-1}}f^p,
\end{equation}
for any $q\in [1,p]$ as can be easily seen, we immediately conclude that the validity of \eqref{eq:1p12} gives equality in \eqref{eq:4p1} in the limit.
For the opposite direction we suppose that we are given $f,F$ such that $0<f^p\leq F$ and $q,p$ for which $1<q<p$. We suppose that we are given a sequence of non negative functions in $L^p$, $(\phi_n)_n$, whose elements satisfy $\int_X\phi_n\,\mr d\mu=f$ and $\int_X\phi^p_n\,\mr d\mu=F$. Then if we define  $A_n=\int_X\phi^q_n\,\mr d\mu$ for every $n\in N$ we may assume, by passing to a subsequence of $(A_n)_n$ that this sequence converges to a fixed constant (note that $(A_n)_n$ is bounded because of the inequality $f^q\leq A_n\leq F^{q/p}$), which we call $A$. By continuity reasons we may also assume that $(A_n)_n$ is constant, that is $\int_X\phi^q_n\,\mr d\mu=A$ for every $n\in N$.

We additionally assume that $(\phi_n)_n$ satisfies equality in the limit in \eqref{eq:4p1} for the choice of $\beta$ which is described above. We now go back to the proof of Theorem \ref{thm:1} and examine where inequalities where used. In these inequalities now we have equality in the limit for our sequence.
The first one that is used is the following

$$ \int_{A(\phi,I)}\phi^q\,\mr d\mu \geq \frac{1}{\alpha_I^{q-1}}\biggl(\int_{A(\phi,I)}\phi\,\mr d\mu\bigg)^q,$$
where $\alpha_I=\mu(A(\phi,I))$ (this is exactly inequality \eqref{eq:3p2}). Additionally the right member of this inequality equals

$$\frac{\left(\int_I\phi\,\mr d\mu - \sum_{J\in S_\phi, J^\star=I}\int_J\phi\,\mr d\mu\right)^q}{\left(\mu(I) - \sum_{J\in S_\phi, J^\star=I}\mu(J)\right)^{q-1}},$$
which in turn is greater or equal than

$$ \mu(I)\frac{y_I^q}{(\tau_I)^{q-1}} - \sum_{\substack{J\in S_\phi\\ J^\star=I}} \mu(J)\frac{y_J^q}{(\beta+1)^{q-1}},$$
where $\tau_I = (\beta+1)-\beta\rho_I$. Since we have equality in the limit in \eqref{eq:4p1}, we conclude that in the inequality

\begin{equation} \label{eq:4p3}
0\leq \sum_{I\in S_\phi}y^{p-q}_I\Biggl\{\int_{A(\phi,I)}\phi^q \,\mr d\mu-\Big[\mu(I)\frac{y_I^q}{(\tau_I)^{q-1}} - \sum_{\substack{J\in S_\phi\\ J^\star=I}} \mu(J)\frac{y_J^q}{(\beta+1)^{q-1}}\Big]\Biggr\}
\end{equation}
we have equality in the limit as $\phi$ moves along $(\phi_n)_n$ and $S_{\phi}$ is replaced by  $S_{\phi_n}$. That is the right member of \eqref{eq:4p3},
tends to zero for our sequence $(\phi_n)_n$.

Additionally, every term on the sum in \eqref{eq:4p3} is non-negative by the comments mentioned right above.
Thus since $y_I\geq f=y_X$ for every $I \in S_{\phi}$, we have that also the following sum tends to zero

\begin{equation} \label{eq:4p4}
0\leq \sum_{I\in S_\phi}f^{p-q}\Biggl\{\int_{A(\phi,I)}\phi^q \,\mr d\mu-\Big[\mu(I)\frac{y_I^q}{(\tau_I)^{q-1}} - \sum_{\substack{J\in S_\phi\\ J^\star=I}} \mu(J)\frac{y_J^q}{(\beta+1)^{q-1}}\Big]\Biggr\},
\end{equation}
as $\phi$ moves along $(\phi_n)_n$. Cancelling the term $f^{p-q}$, using Lemma 2.1 iii) for $I=X$ and the integral assumptions for every $\phi\in(\phi_n)_n$
we immediately conclude that the following inequality is true

\begin{equation} \label{eq:4p5}
A\geq\sum_{I\in S_\phi}\Big[\mu(I)\frac{y_I^q}{(\tau_I)^{q-1}} - \sum_{\substack{J\in S_\phi\\ J^\star=I}} \mu(J)\frac{y_J^q}{(\beta+1)^{q-1}}\Big],
\end{equation}
and is also equality in the limit for our sequence $(\phi_n)_n$. We substitute $\tau_I$ with its value and we get the inequality

\begin{equation} \label{eq:4p6}
A\geq\sum_{I\in S_\phi}\frac{\alpha_I}{\rho_I}\frac{y_I^q}{[(\beta+1)-\beta\rho_I]^{q-1}} -\sum_{I\in S_\phi}\sum_{\substack{J\in S_\phi\\ J^\star=I}} \mu(J)\frac{y_J^q}{(\beta+1)^{q-1}},
\end{equation}
with equality in the limit. Now the right hand side of this inequality equals

$$\sum_{I\in S_\phi}\frac{\alpha_I}{\rho_I}\frac{y_I^q}{[(\beta+1)-\beta\rho_I]^{q-1}} -\sum_{I\in S_\phi,\\ I\neq X}\mu(I)\frac{y_I^q}{(\beta+1)^{q-1}},$$
which in turn equals to

$$\sum_{I\in S_\phi}\frac{\alpha_I}{\rho_I}\Big\{\frac{1}{[(\beta+1)-\beta\rho_I]^{q-1}}-\frac{1}{(\beta+1)^{q-1}}\Big\}{y_I^q}+\frac{y_X^q}{(\beta+1)^{q-1}}.$$
But in the proof of Theorem \ref{thm:1} we have used the inequality
\[
\frac{1}{((\beta+1)-\beta \rho_I)^{q-1}} - \frac{1}{(\beta+1)^{q-1}} \geq \frac{(q-1)\beta \rho_I}{(\beta+1)^q},
\]
for every $I\in S_{\phi}$ and by the same arguments that were used above (by replacing $y_I^p$ by $f^{p-q}y_I^q$) we conclude that we should have equality in the limit in the following inequality

\begin{equation} \label{eq:4p7}
A\geq\sum_{I\in S_\phi}\frac{\alpha_I}{\rho_I}\frac{(q-1)\beta \rho_I}{(\beta+1)^q} y_I^q+\frac{f^{q}}{(\beta+1)^{q-1}}=\frac{(q-1)\beta}{(\beta+1)^q}
\int_X(\mt\phi)^q\,\mr d\mu+\frac{f^{q}}{(\beta+1)^{q-1}}.
\end{equation}
This gives us equality in the limit in the following inequality

\begin{equation} \label{eq:4p8}
\int_X(\mt\phi)^q\,\mr d\mu\leq (1+\frac{1}{\beta})\frac{(\beta+1)^{q-1}A-f^q}{q-1}
\end{equation}
where $\beta$ satisfies $\beta+1=\omega_p(\frac{f^p}{F})$.

But the right side of \eqref{eq:4p8} is minimized exactly when $\beta+1=\omega_q(\frac{f^q}{A})$ as can be seen in \cite{4}, or by a simple calculus argument. From the above we conclude that the value of $A$ satisfies $\beta+1=\omega_q(\frac{f^q}{A})=\omega_p(\frac{f^p}{F})$ and replacing $\beta+1$ by its value in \eqref{eq:4p8} we easily see that

\begin{equation} \label{eq:4p9}
\underset{n}{\lim}\int_X(\mt\phi_n)^q\,\mr d\mu=\omega_q(\frac{f^q}{A})^q A,
\end{equation}
that is, $(\phi_n)_n$ behaves as an extremal sequence for the Bellman function $B_{\mc T}^{(q)}(f,A)$. By using the results of \cite{9} we get that all such sequences behave like $L^q-$ approximate eigenfunctions for the eigenvalue $\omega_q(\frac{f^q}{A})$ which equals $\beta+1$. That is the following holds

\begin{equation} \label{eq:4p10}
\underset{n}{\lim}\int_X|(\mt\phi_n)-(\beta+1)\phi_n|^q\,\mr d\mu=0.
\end{equation}

Our purpose was to show the same equality but with $p$ in place of $q$. This is now not difficult to show because of the following arguments.
Since \eqref{eq:4p10} is true, by a well known theorem in measure theory, we conclude that there exists a subsequence of $(\phi_n)_n$ (without loss of generality we call it again $(\phi_n)_n$) for which $(\mt\phi_n)-(\beta+1)\phi_n\rightarrow 0$ almost uniformly, that is there exists
a decreasing sequence $(A_n)_n$ of $\mu-$ measurable subsets of $X$ for which $\mu(A_n)\rightarrow 0$ and

\begin{equation} \label{eq:4p11}
|(\mt\phi_n)(x)-(\beta+1)\phi_n(x)|\leq \frac{1}{n},
\end{equation}
for every $x\in X\setminus A_n$ and for every $n\in N$. Define now $h_n(x)=(\mt\phi_n)(x)-(\beta+1)\phi_n(x)$ for every $x\in X$.
Then for every $n\in N$
$$\int_X|h_n|^p\,\mr d\mu=\int_{X\setminus A_n}|h_n|^{q}|h_n|^{p-q}\,\mr d\mu+\int_{A_n}|h_n|^p\,\mr d\mu.$$

The first integral of the right side of this last equation is less or equal than $\frac{1}{n^{p-q}}\int_{X\setminus A_n}|h_n|^{q}$, which obviously tends to zero. We proceed now to prove that $\underset{n}{\lim}\int_{A_n}|h_n|^p\,\mr d\mu=0$. By the definition of $h_n$ and since $\mt\phi\geq\phi$ almost everywhere, for every integrable $\phi$ (the tree $\mc T$ differentiates $L^1(X,\mu)$), we see that it is enough to show that
$\underset{n}{\lim}\int_{A_n}(\mt\phi_n)^p\,\mr d\mu=0$.

We define $g_n=\phi^\star_n$, so by Theorem 2.1 we see that
\begin{equation} \label{eq:4p12}
\int_{A_n}(\mt\phi_n)^p\,\mr d\mu \leq \int_0^{\delta_n}\left(\frac{1}{t}\int_0^t g_n\right)^p\mr dt,
\end{equation}
where $\delta_n=\mu(A_n)$ for every $n\in N$. Additionally $g_n$ is equimeasurable with $\phi_n$ so that the $1,q,p$ norms of these two functions are identical, for each $n\in N$. Thus again by Theorem 2.1
\begin{equation} \label{eq:4p13}
\int_{X}(\mt\phi_n)^q\,\mr d\mu \leq \int_0^{1}\left(\frac{1}{t}\int_0^t g_n\right)^q\mr dt \leq A\omega_q(\frac{f^q}{A})^q,
\end{equation}
But by \eqref{eq:4p10} we have that
$\underset{n}{\lim} \int_{X}(\mt\phi_n)^p\,\mr d\mu=A\omega_q(\frac{f^q}{A})^q$, since $\beta+1=A\omega_q(\frac{f^q}{A})^q$.
Thus by \eqref{eq:4p13} we get $\underset{n}{\lim} \int_0^1\left(\frac{1}{t}\int_0^t g_n\right)^q\mr dt=A\omega_q(\frac{f^q}{A})^q$, which means that $(g_n)_n$
is extremal for the respective to the Bellman function problem related to the Hardy operator for the variables $f, A (q>1)$. This gives us in view of the results in \cite{12}, that $(g_n)_n$ tends in the $L^q$-norm to the function $g$, which is defined by $g(t)=\frac{f}{\alpha}t^{-1+\frac{1}{\alpha}}$, $t\in(0,1]$,
where $\alpha=\omega_q(\frac{f^q}{A})^q=\beta+1$. But since $\beta+1=\omega_p(\frac{f^p}{F})^p$, it is easy to see that
$\int_0^1g^p=F$. Obviously
$$ \int_0^{\delta_n}\left(\frac{1}{t}\int_0^t g_n\right)^p\mr dt \leq \left(\frac{p}{p-1}\right)^p\int_0^{\delta_n} g_n^p,$$
so because of \eqref{eq:4p12} it is enough to show that $\lim_{n\to\infty} \int_0^{\delta_n} g_n^p=0$.

It is a standard fact now from measure theory that if for a sequence of integrable functions $(k_n)_n$, defined in a measure space $(Y,r)$, we have that for some integrable function $k$, $\underset{n}{\lim}\int_Y|k_n|\mr dr=\int_Y|k|\mr dr$ and that $k_n$ tends $r$-almost everywhere to $k$ on $Y$, then the sequence $(k_n)_n$ tends to $k$ in the $L^1$-norm (see for example \cite{4}, Theorem 13.47, page 208). Now $\int_0^1g_n^p=\int_0^1g^p=F$ for every $n\in N$ and since $(g_n)_n$ converges in the $L^q$-norm to $g$ we can assume (by passing if necessary to a subsequence) that $g_n$ tends almost everywhere to $g$, so because of the fact mentioned just before we have as a consequence that $g_n^p$ tends to $g^p$ in the $L^1$-norm thus giving us the convergence of $g_n$ to $g$ in the $L^p$-norm, in view of the elementary inequality $(x-y)^p\leq x^p-y^p$ which is true whenever $0\leq y\leq x$, and $p>1$.

Moreover in a finite measure space $(Y,r)$, if we are given a sequence of $p$-integrable functions $(k_n)_n$ and a $p$-integrable function $k$ for which $k_n$ tends almost everywhere and in the $L^p$-norm to $k$ then the following is true
$$\lim_{r(E)\rightarrow 0}\int_E|k_n|\,\mr dr=0,$$
uniformly in $n\in N$. This result is known in the literature as Vitali's convergence theorem and can be seen in \cite{4} (Exercise 13.38, page 203). From all the above we conclude immediately that $\underset{n}{\lim}\int_0^{\delta_n}g_n^p=0$, since $\delta_n=\mu(A_n)\rightarrow0$. Our proof is complete.

\end{proof}

Nikolidakis Eleftherios, Assistant Professor, University of Ioannina, Department of Mathematics, GR 45110, Panepistimioupolis, Greece.

\end{document}